\documentclass[11pt, reqno]{amsart}
\usepackage{geometry,hyperref,amsrefs,  enumitem, cancel}                
\usepackage{color}
\usepackage[normalem]{ulem}

\font\smallsmc = cmcsc9
\font\smalltt = cmtt8
\font\smallit = cmti8

\usepackage{graphicx}
\usepackage{amssymb}
\usepackage{mathbbol}
\usepackage{epstopdf}
\usepackage{verbatim}
\usepackage[mathscr]{euscript}

\numberwithin{equation}{section}

\theoremstyle{plain}
\newtheorem{theo}{Theorem}[section]
\newtheorem{lem}[theo]{Lemma}
\newtheorem{prop}[theo]{Proposition}

\newtheorem{cor}[theo]{Corollary}

\theoremstyle{definition}

\newtheorem{definition}[theo]{Definition}

\theoremstyle{plain}

\theoremstyle{definition}

\newcommand{\beq}{\begin{equation}}
\newcommand{\eeq}{\end{equation}}
\newcommand{\beqn}{\begin{equation*}}
\newcommand{\eeqn}{\end{equation*}}
\renewcommand{\a}{\alpha}
\renewcommand{\b}{\beta}

\renewcommand{\d}{\delta}

\newcommand{\f}{\varphi}
\newcommand{\g}{\gamma}
\newcommand{\h}{\eta}

\renewcommand{\o}{\omega}
\renewcommand{\q}{\vartheta}

\newcommand{\s}{\sigma}

\newcommand{\G}{\Gamma}

\newcommand{\bR}{\mathbb{R}}

\newcommand{\cH}{\mathscr{H}}

\newcommand{\cK}{\mathscr{K}}
\newcommand{\cL}{\mathscr{L}}
\newcommand{\cM}{\mathscr{M}}

\newcommand{\cQ}{\mathscr{Q}}

\newcommand{\cU}{\mathscr{U}}
\newcommand{\cV}{\mathscr{V}}
\newcommand{\cW}{\mathscr{W}}

\renewcommand{\square}{\kern1pt\vbox
{\hrule height 0.6pt\hbox{\vrule width 0.6pt\hskip 3pt
\vbox{\vskip 6pt}\hskip 3pt\vrule width 0.6pt}\hrule height0.6pt}\kern1pt}

\renewcommand\={:=}

\newcommand{\wt}{\widetilde}
\newcommand{\wh}{\widehat}

\newcommand{\bt}{\begin{theo}\ \ }
\newcommand{\et}{\end{theo}}
\newcommand{\bp}{\begin{prop}\ \ }
\newcommand{\ep}{\end{prop}}
\newcommand{\bc}{\begin{cor}\ \ }
\newcommand{\ec}{\end{cor}}
\newcommand{\bl}{\begin{lem}\ \ }
\newcommand{\el}{\end{lem}}
\newcommand{\bd}{\begin{definition}}
\newcommand{\ed}{\end{definition}}
\newcommand{\n}{\nabla}

\newcommand{\be}{\begin{equation}}
\newcommand{\ee}{\end{equation}}

\def\<#1,#2>{\langle\,#1,\,#2\,\rangle}
\newcommand{\arr}{\begin{array}{rlll}}
\newcommand{\ea}{\end{array}}
\newcommand{\bea}{\begin{eqnarray}}
\newcommand{\eea}{\end{eqnarray}}
\newcommand{\bean}{\begin{eqnarray*}}
\newcommand{\eean}{\end{eqnarray*}}

\newcommand{\Ga}[3]{\G_{#1 #2}^{\phantom{#1} #3}}
\newcommand{\GGa}[3]{{\mathbf \G}_{#1 #2}^{\phantom{#1} #3}}

\def\sideremark#1{\ifvmode\leavevmode\fi\vadjust{
\vbox to0pt{\hbox to 0pt{\hskip\hsize\hskip1em
\vbox{\hsize3cm\tiny\raggedright\pretolerance10000
\noindent #1\hfill}\hss}\vbox to8pt{\vfil}\vss}}}

\renewcommand{\sf}{shearfree }

\newcommand{\ps}{{\operatorname{p}}}
\newcommand{\qs}{{\operatorname{q}}}
\newcommand{\grad}{{\operatorname{grad}}}

\title[The Levi-Civita connections of  manifolds with prescribed optical geometries]{The Levi-Civita connections of  Lorentzian manifolds \\ with prescribed optical geometries}
\author{Dmitri V. Alekseevsky, Masoud Ganji, Gerd Schmalz and Andrea Spiro}

\begin{document}
\begin{abstract} 
We explicitly derive the Christoffel symbols in terms of adapted frame fields for the Levi-Civita connection of a Lorentzian $n$-manifold $(M, g)$, equipped with a prescribed optical geometry of  K\"ahler-Sasaki  type. The formulas found in this paper have several important applications, such as determining the geometric invariants of Lorentzian
manifolds with prescribed optical geometries or  solving curvature constraints.
\end{abstract}
\thanks{D. V. Alekseevsky was supported by the Grant "Basis-foundation
(Leader)" 22-7-1-34-1.
}
\maketitle
\section{Introduction}
An  {\it optical geometry}, a notion  introduced in  the late eighties by Robinson and Trautman,  is a geometrical structure  that encodes 
 the existence of an electromagnetic plane wave -- or  an appropriate higher dimensional generalisation \cite{AGSS} --  propagating along a prescribed foliation  by curves of a Lorentzian  manifold. 
 Let us recall the relevant definitions.  A  {\it null congruence} on a Lorentzian $n$-manifold  $(M, g)$, $n \geq 3$,  is  a   foliation   by    curves, which are  tangent to  some nowhere vanishing   null  vector  field. 
  Given a Lorentzian $n$-manifold  $(M, g)$, $n \geq 3$, 
  a   null congruence  is called {\it geodesic shearfree},  or  {\it shearfree} for short,  if there is a  choice for  a  nowhere vanishing  tangent null vector field $\ps$,   whose local flow  preserves both the  codimension  one   distribution $ \cW:=\ps^{\perp_g}$  and the  conformal class of the induced  degenerate metric $ h \= g|_{\cW \times \cW}$  on the spaces $\cW_x = \ps^{\perp_g}|_x$, $x \in M$.   These conditions are  equivalent to  requiring  that the Lie derivative $\cL_\ps g$  has the form
\beq \label{defsf2}   \cL_{\ps} g = f g + {\ps}^\flat \vee \h\qquad \text{for some function $f$ and some $1$-form}\ \ \h\ .\eeq
If this holds, the vector field  $\ps$ is also geodesic, i.e. $\nabla_{\ps}\ps =\lambda \ps$, and the curves of the congruence are geodesics (see e.g.  \cites{RT, AGS, AGSS, FLT}).  A quadruple $\cQ \= (\ps, \cW, [h], \{g\})$,  given by
\begin{itemize}[leftmargin = 20 pt]
\item[(a)]  a nowhere vanishing  vector field $\ps$, determined up to multiplication by a nowhere smooth function $f$, 
\item[(b)] a codimension one distribution $\cW$, 
\item[(c)] a conformal class of semi-positive  metrics on $\cW$, 
\item[(d)] a non-empty set of  Lorentzian metrics $\{g\}$,  which are exactly all   metrics $g$  with respect to  which $\ps$ is a  null  vector field with  $\cW = \ps^{\perp_g}$ and $[h] = [g|_{\cW \times \cW}$ and both $\cW$ and $[h]$  are preserved by the local flow of $\ps$. 
\end{itemize}
is an  {\it optical geometry} in the sense of Robinson and Trautman \cites{RT, AGSS, FLT} 
 (\footnote{As a matter of fact,  all four  elements of   $\cQ$   can be   recovered just by (i)  the  $1$-dimensional distribution $\cK$, which is generated by  $\ps$  and  (ii) the set of metrics $\{g\}$, provided that they satisfy appropriate conditions.  Thus, the optical geometries can be also   defined  as such  pairs  $(\cK, \{g\})$ -- see the original definition  in   \cite{RT}.}).  The Lorentzian metrics  $g$ in  the  set $\{g\}$   are called {\it  compatible with  the prescribed  optical geometry $\cQ$}.  \par
\smallskip
By Robinson's Theorem \cites{Ro, HM},  the  \sf  null congruences of  a real analytic four dimensional  Lorentzian manifold are exactly  the foliations by the lines of  propagation of  electromagnetic plane waves. \par
\smallskip
Many interesting examples of  optical geometries $\cQ = (\ps, \cW, [h], \{g\})$ are provided by  connections on principal  $A$-bundles  $\pi: M\to S = M/A$  with one-dimensional structure groups $A = \bR$ or $S^1$. On each bundle of this kind, one may consider an optical geometry in which $\ps$ is  the generator of  the action of the group $A$ along the fibers,  and $\cW$ and  $[h]$  are the appropriate $A$-invariant  distribution and conformal class. In this case, the quadruple  $\cM \= (\pi: M \to S, \ps, \cW, [h])$ is called a {\it  regular \sf manifold}  
and a metric $g \in \{g\}$ of the corresponding optical geometry $\cQ = (\ps, \cW, [h], \{g\})$ is said to be a  {\it compatible metric} of $\cM$. \par
\smallskip
The regular \sf manifolds are important geometric objects not only for their role in Lorentzian geometry, but also for their relations with CR geometry. Indeed,  for any  regular \sf manifold $\cM \= (\pi: M \to S, \ps, \cW, [h])$, 
 the  base manifold  $S = M/A$  is naturally 
equipped with  the     codimension one distribution $\cW^S \subset TS$ and the  positive definite conformal  metric $[h^S]$ that are  obtained by  projecting the  $A$-invariant distribution  $ \cW:=\ps^{\perp}$  and conformal class $[h]$ onto  $S = M/A$.   If  $M$ is even dimensional and the projected   
   distribution  $\cW^S \subset TS$ is  contact then the regular \sf manifold $\cM$ is called {\it (maximally) twisting}. For any  such $\cM$,  
 the corresponding optical geometry  $\cQ = (\ps, \cW, [h], \{g\})$   determines a  family  $J^S$  of complex structures $J^S_x: \cW^S_x \to \cW^S_x$  on the projected distribution of $S$,  that make $S$ a  {\it strongly pseudoconvex almost CR manifold}  (see, e.g. \cites{HLN, AGS, AGSS, FLT} and references therein). \par
 \smallskip
Celebrated examples of twisting  regular  \sf manifolds  are given by the  $4$ dimensional  space-times with  Taub-NUT metrics and the $4$-dimensional   Kerr black holes.  For such Lorentzian manifolds,  
the base manifold  of the $A$-bundle $\pi: M \to S$ has an additional remarkable geometric feature: it is a principal bundle $\pi^S = S \to N$ with one dimensional structure group $A' = \bR$ or $A' = S^1$, and the base manifold  $N = S/A' = M/(A \cdot A')$ has a natural structure of a
  K\"ahler manifold. Moreover,  the strongly pseudoconvex almost CR manifold $(S, \cW^S, J^S)$ is   a {\it regular Sasaki manifold}  and  the structure  group $A'$   of  $S$ preserves
\begin{itemize}[leftmargin = 20 pt]
\item[(i)] the CR structure $(\cW^S, J^S)$, 
\item[(ii)] a contact $1$-form $\theta_o$ for $\cW^S$, i.e.,     $\cW^S = \ker \theta_o$,  such that  $d \theta_o = \pi^{S*} \o_o$ for some K\"ahler form  $\o_o = g_o(J \cdot, \cdot)$  on  $(N, J)$; 
 \item[(iii)] the conformal class $[h]$ on $\cW$ contains the degenerate metric $h_o = \big((\pi^S \circ \pi)^* g_o\big)|_{\cW}$. 
\end{itemize}
The fact that the Taub-NUT  and Kerr metrics have these properties is one of the reasons of the  interest in  twisting regular \sf  manifolds,  in which the  almost  CR manifold $(S, \cW^S, J^S)$ is a Sasaki manifold   projecting  onto a K\"ahler manifold.  Such manifolds  are called  {\it of K\"ahler-Sasaki type}   \cite{AGSS}. \par
 \medskip 

As is well known, the geometric properties of a Lorentzian metric  $g$ are fully determined by its Levi-Civita connection $\nabla^g$, which in turn  can be expressed through its Christoffel symbols in some fixed frame field. Once the Christoffel symbols are known, one can reconstruct all components of the Riemann curvature tensor $R$ and of  its covariant derivatives  $\n^g \n^g \ldots \n^g R$ in the considered frame field. In principle, these components determine all the geometric properties of $g$. This is because, according to classical results in the theory of $G$-structures, any isometric invariant of a pseudo-Riemannian manifold depends on the components in orthonormal bases of the Riemann  tensor and its covariant derivatives up to an appropriate order. For more information, see,  for example, \cites{St, PS, Ke, Ke1,KS} and references therein.

This observation indicates that the explicit expressions of the Christoffel symbols in appropriate frame fields represent a  fundamental  tool for   studying  the compatible  metrics of a given regular \sf manifold of K\"ahler-Sasaki type  and possibly finding   solutions of the Einstein (or other physically relevant) equations   in this  class of metrics. \par
\smallskip
 In this paper,  we discuss in great detail  the Christoffel symbols of the Levi-Civita connection $\n^g$ of a  compatible metric $g$  of a  regular  \sf manifold $\cM \= (\pi: M \to S, \ps, \cW, [h])$ of K\"ahler-Sasaki type. More precisely, we  
fix a special (locally defined) frame field $(e_1, \ldots, e_n)$, which is   well adapted to the optical geometry and  is    determined  only up to  a choice of a local frame field  on the underlying K\"ahler manifold $N = M/(A \cdot A')$.  Such a frame field  has the following  two  useful properties: 
    \begin{itemize}[leftmargin = 20pt]
\item[(1)] the last two vector fields $e_{n-1}$, $e_{n}$ are the  generators of the actions of the groups $A$ and $A'$, respectively,  and are therefore canonically associated with the considered manifold; 
\item[(2)] the   vector fields $e_i$, $1 \leq i \leq n-2$,  are  tangent to the distribution $\cW$ at all points and are $A\cdot A'$-invariant, thus projecting onto a frame field $(\wt e_1, \ldots, \wt e_{n-2})$ on $N$. 
\end{itemize}
 Note that (1) and (2) allow  to minimise the number of parameters that are necessary to determine  the components of a compatible metric $g$.   Notice also that, due to the fact that $\cM$ is twisting,  a frame field  satisfying  (1) and (2) cannot coincide with a coordinate frame field.  This forces us to avoid  the use of  coordinates in all subsequent computations. \par
 \smallskip
    After choosing an adapted frame field of this kind, we write down the general expression of a  compatible metric $g$  in terms of  its   dual frame field and we  determine  the Christoffel symbols of  $\n^g$ in such  frame and coframe fields,  using just Koszul's formula and  classical results on transformations of Levi-Civita connections under conformal transformations. \par
\smallskip
The expressions for the  Christoffel symbols given in this paper have been originally determined during the preparation of \cite{AGSS} 
and have been successfully  used to derive a coordinate-free characterisation of  the generalised Taub-NUT  metrics on  even dimensional manifolds (see  e.g. \cite{AC} and references therein  for other characterisations of the metrics of such a kind). However,  the details of the  actual computations did not appear in \cite{AGSS} and some formulas of that paper  had  some minor sign errors -- very few indeed and with  no effect on any statement  and proof.   
The same explicit (and amended) expressions  have been later used in \cite{GGSS} for determining explicit expressions for  the  components of the  Ricci tensor  of  compatible  metrics of a  \sf manifold $\cM$ of K\"ahler-Sasaki type satisfying  conditions that  generalise Kerr's ansatz  for the  $4$-dimensional rotating black holes.   These  expressions for the Ricci tensor allowed us to  translate  the Einstein equations for a compatible metric into equations  on  its parameters in an adapted frame and to find a large  class  of exact solutions  that  naturally includes the classical  Kerr black holes.  
 We anticipate a number of further applications of the explicit expressions of these Christoffel symbols and believe that  the detailed computations we present in this paper will be a helpful tool for other researchers who are interested in the developments of this  field.
 
  \par
\smallskip
%
%
The paper is structured into two sections:  In \S 2, we define  the adapted frame fields of a compatible metric,  that is the frame fields  in which all computations of this paper are performed; In  \S  3 we derive the explicit  list of Christoffel symbols  and provide the details of the computations.\par
\medskip
\section{The general form of a compatible metric\\ on a \sf manifold of K\"ahler-Sasaki type}

\subsection{Notational issues}
Consider a    \sf manifold $\cM \= (\pi: M \to S,$ $\ps, \cW,$ $[h])$ of K\"ahler-Sasaki type.  We  use the following notation: 
\begin{itemize} [leftmargin = 20pt]
\item[(1)]  $(N, J, g_o)$ is the K\"ahler manifold onto which $S$ projects and $\o_o = g_o(J \cdot, \cdot)$  is  the K\"ahler form of $N$ (\footnote{Note that there is  a sign difference in the  definition of  $\o_o$   w.r.t. \cite{AGSS}. There   it is  defined as
 $\o_o := g_o(\cdot , J\cdot)$.}); 
\item[(2)] $A$ and $A'$ are the $1$-dimensional structure groups of the principal bundles $\pi: M \to S$ and $\pi^S: S \to N$, respectively; 
\item[(3)]  $\ps_o$   and $\qs_o^S$ are  the   fundamental  vector fields of the principal bundles  $\pi: M \to S$  and $\pi^S: S \to N$,  corresponding to the element of the standard basis of $Lie(A) = Lie(A') = \bR$.  This means that   $\Phi^{\ps_o}_s(x) = e^s(x)$, $x \in M$,  and   $\Phi^{\qs^S_o}_s(y) = e^s(y)$, $y \in S$;  
\item[(4)] $\theta_o$ is the contact $A'$-invariant $1$-form on $S$ satisfying the conditions
\beq d \theta_o = \pi^{S*} \o_o \ , \qquad \theta_o(\qs_o) = 1\ , \qquad  \ker \theta_o|_x = \cW^S_x\  \ , \ x \in S\ ; \eeq
and $\q_o$ is the pull-back $ \q_o = \pi^*(\theta_o)$ of $\theta_o$ on $M$.  
\end{itemize}
\par
It is important to note that   $\cW^S$ is an $A'$-invariant  horizontal distribution on the principal bundle $\pi^S : S \to N$, and it is therefore  a connection for this bundle. The associated  connection $1$-form is   $\theta_o$ and  its curvature $2$-form is  $d \theta_o =  \pi^{S*} \o_o$.  \par
\smallskip
For what concerns  the $A$-bundle $\pi: M \to S$, throughout the paper  {\it we assume that it  is  trivial and equipped with the natural  {\rm flat} connection of a Cartesian product}.  This apparently restrictive condition can be always locally satisfied replacing $S$ by  an open subset $\cV \subset S$,  on which the bundle is trivialisable,   and identifying $\pi: M \to S$ with  the trivial bundle $\pi: \pi^{-1}(\cV) \simeq \cV \times A \to \cV$ equipped with the standard flat connection. \par
\smallskip
We denote by $\cH_o$ the horizontal distribution of the flat connection of $\pi: M \to S$.\par
\smallskip
For any given vector field $X$ on the K\"ahler manifold $N$,  we denote by
\begin{itemize}[leftmargin = 20pt]
\item[--] $X^{(S)}$ the unique $A'$-invariant  horizontal vector field in $\cW^S \subset S$ projecting onto $X$; 
\item[--] $\wh X$ the  unique $A$-invariant horizontal  vector field in $\cH$ projecting onto $X^{(S)}$ and thus also onto  $X$; note that,  by  definition of $\cW^S$,  the vector field $\wh X$  takes 
values in $\cH_o \cap \cW$.
\end{itemize}
The unique  $A$-invariant horizontal vector field in $\cH_o$ projecting onto   $\qs_o^S$  is  denoted by $\qs_o$.
\par
\smallskip
Owing to the  $A$- and $A'$- invariance of the connections of $\pi: M \to S$ and $\pi^S: S \to N$ and  the properties of the connection $1$-form $\theta_o$,  for any pair of  vector fields $X, Y $ on $N$ the following Lie bracket relations hold (\footnote{The Lie bracket $[\wh X, \wh Y]$ differs  by a sign  from the one used  in  \cite{AGSS}. Since in both papers, it is assumed   $d \theta_o = \o_o$, the sign difference is a consequence of the  different definitions of the K\"ahler form $\o_o$.}):
\beq \label{62}  [\wh X, \wh Y]  - \wh{[X, Y]} =  - g_o(JX, Y) \qs_o\ ,\qquad [\wh X, \ps_o] = [\wh X, \qs_o] =  [\ps_o, \qs_o] = 0\ .
\eeq
%
%

\smallskip
\subsection{The adapted frame fields}
Consider a frame field $(E_1, \ldots, E_{n-2})$ on an open set $\cV \subset N$ of the K\"ahler manifold and the corresponding lifted vector fields $(\wh E_1, \ldots, \wh E_{n-2})$ on $M$, taking values in the distribution $\cW' = \cH \cap \cW$. The vector fields of the $(n-1)$-tuple $(\wh E_1, \ldots, \wh E_{n-2}, \ps_o)$ are pointwise linearly independent and hence give  linear frames for  the  spaces $\cW_x  \subset T_x M$, $x \in \cU = (\pi^S \circ \pi)^{-1}(\cV)$. Since $\qs_o$ projects onto $\qs_o^S$ and $\qs_o^S$ is transversal to $\cW^S = \pi_*(\cW)$,   the vector fields of the $n$-tuple 
\beq \label{adapted} (\wh E_1, \ldots, \wh E_{n-2}, \ps_o, \qs_o)\eeq
are  pointwise linearly independent and  determine  a frame field on   $ \cU$.  We call \eqref{adapted} the {\it adapted frame field of $\cM$ determined by the  frame field $(E_i)$ on $N$}.\par
Note that, due to \eqref{62}, the Lie brackets between any two vector fields of an adapted frame have the form
\beq \label{62*}  [\wh E_i, \wh E_j]  = c_{ij}^k \wh E_k  - g_o(JE_i, E_j) \qs_o\ ,\qquad [\wh E_i, \ps_o] = [\wh E_i, \qs_o] =  [\ps_o, \qs_o] = 0\ ,
\eeq
where  the  $c_{ij}^k$ are the   functions  such that 
$[E_i, E_j] = c_{ij}^k E_k$.\par
The dual coframe field of  $(\wh E_1,\ldots, \wh E_{n-2},  \ps_o, \qs_o)$
is denoted by $(\wh E^1, \ldots, \wh E^{n-2}, \ps_o^*, \qs_o^*)$.  Since the dual $1$-form $\qs_o^*$ satisfies  $\qs_o^*(\qs_o) = 1$ and vanishes identically on $\cW$ (because  $\cW$ is spanned by the $\wh E_i$ and $\ps_o$),  it has the same kernel and takes the same value on $\qs_o$ as the $1$-form $\q_o$. Thus 
\beq \label{obs} \qs_o^* = \q_o\eeq
{\it for {\rm  any} choice of the adapted frame $(\wh E_i, \ps_o, \qs_o)$}.\par
\medskip

\subsection{Parameterisation of the compatible  metrics} Let $(E_i)$ be a (local) frame field on $N$ and denote by $(\wh E_1,\ldots, \wh E_{n-2},  \ps_o, \qs_o)$ the corresponding adapted frame field for $\cM$.  Since we are assuming that $\cM$ is of K\"ahler-Sasaki type,  the conformal class $[h]$   consists of  the degenerate metrics on $\cW$ having the form
\beq \label{metric1} h = \s (\pi^S \circ \pi)^*(g_o)|_{\cW} ,\qquad \s = \text{conformal scaling factor}\ .\eeq
By the results in \cite{AGSS}*{\S 2.5} (see also \cite{GGSS}),  the compatible  Lorentzian  metrics  on $\cM$ are locally  in one-to-one correspondence with the pairs $(h, \qs)$ given by 
\begin{itemize}[leftmargin = 22 pt]
\item  a degenerated metric $h$  on $\cW$ as in    \eqref{metric1}: 
\item a vector field $\qs$, which is transversal to the distribution $\cW = \cW' + \bR \ps_o$, i.e., of the form  
\beq \label{52bis} \qs \=a \qs_o + b \ps_o + c^i \wh E_i\ ,\qquad a \neq 0\ .\eeq
\end{itemize}
More precisely, given  the conformal factor $\s$ and the vector field $\qs$, the  corresponding  compatible metric $g = g^{(\s, \qs)}$ is  the unique Lorentzian metric satisfying  conditions
\beq  \label{52}
\begin{split}
&g(\wh X, \wh Y) = \s g_o(X, Y)\ ,\qquad g(\wh X, \ps_o) = g(\ps_o, \ps_o) = 0\ , \\[10pt]
 &g(\wh X, \qs)  =  0\ ,\qquad g(\ps_o, \qs) = 1\ ,\qquad g(\qs, \qs)   = 0\ .
 \end{split}
 \eeq
From \eqref{52bis} and the first line of \eqref{52},  the second line  in \eqref{52}  is  equivalent to
\beq
\begin{split}
 &g(\wh X, \qs_o) = - \frac{c^i \s}{a} g_o(X, E_i)\ ,\qquad  g(\ps_o, \qs_o)  = \frac{1}{a}\ ,\\
 &g( \qs_o ,  \qs_o ) = - 2 \frac{b}{a^2 }  +  \frac{1}{a^2}  c^i c^j  \s g_o( E_j , E_i) \ .
\end{split}
\eeq
 Introducing the shorter notation  
\beq \label{solve}  \a \= \frac{2}{a \s} \ ,\qquad \b \=  \frac{2}{\s}\left( - 2 \frac{b}{a^2} +  \frac{1}{a^2} c^i c^j \s g_o( E_j , E_i)\right)  \ ,\qquad  \g^i \= -  2  \frac{c^i}{a} \ ,\eeq
we get  that   $g = g^{(\s, \qs)}$ is  the unique Lorentzian metric satisfying the condition
 \beq \label{condiriso}
\begin{split}
& g(\wh X, \wh Y) = \s g_o(X, Y)\ ,\qquad g(\wh X, \ps_o) =  g(\ps_o, \ps_o) =  0\ ,\qquad g(\ps_o, \qs_o) = \frac{\s\a}{2}\ ,\\[5pt]
&  g(\qs_o, \wh X) = \frac{\s \g^i}{2}  g_o(X, E_i)\ ,\qquad g( \qs_o ,  \qs_o ) = \frac{\s}{2} \b\ .
\end{split}
\eeq
This means that   $g$ has the form
\beq \label{buona}
\begin{split} g & =  \s g_o(E_i, E_j)  \wh E^i \vee \wh E^j+ \\
&\hskip  3 cm  +\qs^*_o \vee \left(\s \a \ps_o^* +  \s  \g^i  g_o(E_i, E_k) \wh E^k  + \frac{\s \b}{2} \qs^*_o\right) = \\
& =   \s \bigg\{( \pi^S \circ \pi)^*(g_o)\big|_{\cW'}+ \\
&\hskip  3 cm  +\q_o \vee \left( \a \ps_o^* +   \g^i g_o(E_k, E_i) \wh E^k  +  \frac{\b}{2}  \q_o\right)\bigg\}.
\end{split}
\eeq
The expression \eqref{buona} gives a convenient parameterisation in terms of the $(n+1)$-tuple of smooth functions $(\s, \a, \b, \g^i)$ for the    compatible  metrics of $\cM  = (\pi: M \to S, \ps, \cW, [h])$. We emphasise that, conversely, any metric having the form  \eqref{buona},  for some $\s > 0$ and $\a \neq 0$, is a compatible  metric. Indeed, it  is  associated with the conformal factor $\s$ and the vector field $\qs = a \qs_o + b \ps_o + c^i \wh E_i$ where  $a$, $b$ and $c^j$ are solutions to \eqref{solve} for the given $\a$, $\b$ and $\g^i$. They are 
$$a = \frac{2}{\a \s}\ ,\qquad  b \= - \frac{ \b}{\a^2 \s} +  \frac{1}{2\a^2 \s} \g^i \g^j  g_o( E_j , E_i)    \ ,\qquad  c^i = - \frac{\g^i}{\a \s}\ .$$
\par 
\medskip
\section{The Christoffel symbols  in an adapted frame field of the Levi-Civita connection of a compatible Lorentzian  metric}
\subsection{The  complete list  of  the Christoffel symbols}
Let  $\cM  = (\pi: M \to S, \ps, \cW, [h])$ be a twisting regular \sf manifold of K\"ahler-Sasaki type, with $S$ projecting onto the K\"ahler manifold $(N, J, g_o)$. 
Let also  $(E_i)$ be a frame field  on an open set $\cV \subset N$ and  $(X_A) = (\wh E_1,\ldots, \wh E_{n-2},  \ps_o, \qs_o)$ the corresponding adapted frame field on $\cU = (\pi^S \circ \pi)^{-1}(\cV) \subset M$.  We use the notation $g_{ij}$, $\o_{ij}$, $J_i^j$,  $c_{ij}^k$ for  the functions defined by 
$$g_{ij} \= g_o(E_i, E_j)\ ,\qquad \o_{ij} \= g_o(J E_i,  E_j)\ ,\qquad J E_i = J_i^j E_j\ ,\qquad [E_i, E_j] = c_{ij}^k E_k\ .$$
For what concerns the Christoffel symbols $\GGa A B C$ (i.e., the functions  defined by  $\n_{X_A} X_B = \GGa A B C X_C$), we are going to use  the convention that    $\GGa i j m$   denotes the function which gives the component  of   $\n_{\wh E_i} \wh E_j$  in the direction of $\wh E_m$,   $\GGa i j {\ps_o}$    is the function that gives    the component of  $\n_{\wh E_i} \wh E_j$  in the direction of $\ps_o$,     $\GGa i j {\qs_o}$   is  the function giving the component of  $\n_{\wh E_i} \wh E_j$ in the direction of $\qs_o$,  and so on.
\par 
\smallskip
Our main result is the following:
\begin{theo} \label{main} Let $g$ be a compatible metric for $\cM$,  hence of the form \eqref{buona} for an  $(n+1)$-tuple of smooth functions  $(\s, \a, \b, \g^i)$ on $\cU$, with $\s > 0$ and $\a \neq 0$ at all points.  The  Christoffel symbols $\GGa A B C$ of the Levi-Civita connection  of $g$ in the frame field $(X_A)  =  (\wh E_i, \ps_o, \qs_o)$  are given by
\begin{align}
\label{746ter}
\GGa i j m &=
  g^{mk} g_o(\n^o_{E_i}   E_j, E_k)  +  g^{mk}S_{ij|k}
  +  \frac{\g^m \o_{ij} }{4}   + \frac{1}{2 \s}  \wh E_i(\s) \d_j^m +  \frac{1}{2 \s}  \wh E_j(\s) \d_i^m\nonumber \\
 &\hskip 1 cm- \frac{g_{ij}}{2 \s}   \left(  g^{mk} \wh E_k( \s) - \frac{\g^m}{\a} \ps_o(\s)\right)\ , \\
 \nonumber & \qquad \text{where} \ \ S_{ij|k}\ \text{is defined by}\\
 \nonumber  S_{ij|k} &\= 
 \frac{\g^\ell}{4}g_o(J E_i, E_k) g_o(E_\ell, E_j) + \frac{\g^\ell}{4} g_o(J E_j,  E_k)  g_o(E_\ell, E_i) -  \frac{\g^\ell}{4} g_o(J E_i, E_j)  g_o(E_\ell,E_k)\ ,
\\[10pt]%
 \GGa i j {\ps_o}  &=  \frac{1}{2 \a} \wh E_i(\g^k g_{jk} )+  \frac{1}{2 \a}\wh E_j(\g^k g_{ik} )   
-  \frac{1}{4 \a} \g^m \g^k  g_{mk} \o_{ij}   - \frac{ \g^m}{\a} g_o(\n^o_{E_i}   E_j, E_m)  - \frac{ \g^m}{\a}  S_{ij|m} \nonumber\\
 &\hskip 1 cm- \frac{g_{ij}}{2 \s}  \left(\frac{2}{\a} \qs_o(\s) +  \frac{1}{\a^2} \left(  \g^m \g^k  g_{mk} - 2\b \right) \ps_o(\s) - \frac{ \g^m}{\a} \wh E_m(\s)\right)\ , \\
 \GGa i j {\qs_o}  &=  - \frac{  \o_{ij}}{2}  - \frac{ g_{ij} }{\a \s} \ps_o(\s)\ , \\[7 pt]
 \GGa i {\ps_o} m & = \GGa  {\ps_o} i m =  \frac{\a g^{mk} \o_{ik} }{4} +\frac{1}{2 \s}  \ps_o(\s) \d_i^m \ , \\
 \GGa i {\ps_o} {\ps_o}&=  \GGa {\ps_o} i {\ps_o} = \frac{1}{2\a}  \wh E_i(\a) +  \frac{1}{2\a}  \ps_o( \g^k) g_{ik}  -   \frac{ \g^m \o_{im}}{4} + \frac{1}{2 \s} \wh E_i(\s)\ ,  
 \end{align}
 \begin{align}
 \GGa i {\ps_o} {\qs_o}&=  \GGa  {\ps_o} i {\qs_o} = 0\ ,\\[7 pt]
 \GGa i {\qs_o} m &=  \GGa  {\qs_o} i m  = \frac{g^{mk} }{4}\wh E_i(\g^tg_{tk})- \frac{g^{mk} }{4}\wh E_k(\g^tg_{ti})   - \frac{\g^\ell}{4} c_{i r}^t g_{t \ell} g^{m r}+  \frac{g^{mk}\o_{ik} } {4} \b - \nonumber \\
& \hskip 0.5 cm -\frac{\g^m}{4\a}\wh E_i(\a) + \frac{\g^m}{4\a}\ps_o(\g^t )g_{ti} + \frac{1}{2 \s} \qs_o(\s)\d_i^m-   \frac{\g^t}{4 \s }  g_{ti} \left(g^{mk} \wh E_k(\s) - \frac{\g^m}{\a} \ps_o(\s)\right)\ ,\\
 %
 \GGa i {\qs_o} {\ps_o} &=  \GGa  {\qs_o} i  {\ps_o}  =   \frac{1}{2\a}\wh E_i(\b)+  \frac{1}{4 \a^2} \g^m \g^k g_{mk} \wh E_i(\a) - \frac{1}{4 \a^2} \g^m \g^k g_{mk}\ps_o(\g^t )g_{it}  - \frac{1}{2 \a^2} \b\wh E_i(\a) + \nonumber   \\
& \hskip 1 cm +  \frac{1}{ 2\a^2} \b \ps_o(\g^t )g_{it} - \frac{\g^m}{4\a} \wh E_i(\g^tg_{tm}) +  \frac{\g^m}{4\a}\wh E_m(\g^tg_{it})  + \frac{\g^m \g^t}{4 \a} g_{t\ell} c^\ell_{im} +  \frac{\g^m}{4\a}\o_{im} \b - \nonumber \\
 & \hskip 1 cm-    \frac{\g^t}{4 \s} g_{ti} \left(\frac{2}{\a} \qs_o(\s) +  \frac{1}{\a^2} \left(  \g^m \g^k  g_{mk} - 2\b \right) \ps_o(\s) - \frac{ \g^m}{\a} \wh E_m(\s)\right)\ , 
 \\
 \GGa i {\qs_o} {\qs_o} &=  \GGa  {\qs_o}  i {\qs_o} =  \frac{1}{2\a}\wh E_i(\a)- \frac{1}{2\a}\ps_o(\g^t)g_{it}  + \frac{1}{2 \s}  \wh E_i(\s)   -     \frac{\g^t  g_{ti}}{2 \a \s}  \ps_o(\s)\ ,\\
 %
 \GGa {\ps_o}  {\ps_o}  m & =  0\ ,\\
 \GGa {\ps_o}   {\ps_o} {\ps_o} & =  \ps_o(\log ( \a \s))\ , \\
 \GGa {\ps_o}   {\ps_o} {\qs_o} & =  0 \ , \\[7pt]
 \GGa {\ps_o}  {\qs_o}  m & =  \GGa {\qs_o}   {\ps_o}  m = \frac{1}{4}\ps_o(\g^m) - \frac{g^{mk}}{4 } \wh E_k(\a)-  \frac{\a }{4\s} \left (g^{mk} \wh E_k(\s) - \frac{\g^m}{\a} \ps_o(\s) \right)\ ,
\\
 \GGa {\ps_o}   {\qs_o} {\ps_o} & =  \GGa  {\qs_o}  {\ps_o}  {\ps_o}  =  \frac{1}{2\a}\ps_o(\b) -\frac{\g^m}{4\a}\ps_o(\g^i)g_{im}+ \frac{\g^m}{4\a}\wh E_m(\a)+\frac{1}{2 \s} \qs_o(\s)-\nonumber\\
 & \hskip 1 cm- \frac{1}{2 \s}\left(\qs_o(\s) +  \frac{1}{ 2 \a} \left(  \g^m \g^k  g_{mk} -2\b \right) \ps_o( \s) - \frac{ \g^m}{2 } \wh E_m(\s) \right) \ , \\
 \GGa {\ps_o}   {\qs_o} {\qs_o} &  =   \GGa  {\qs_o}  {\ps_o}   {\qs_o} = 0 \ ,
\\[7pt]
\GGa {\qs_o}  {\qs_o}  m & = \frac{g^{mk}}{2}\qs_o(\g^i)g_{ik} - \frac{g^{mk}}{4}\wh E_k(\b)  -\frac{\g^m}{2\a} \qs_o(\a) +\frac{\g^m}{4\a}\ps_o(\b)-\nonumber\\
& \hskip 1 cm-  \frac{\b }{4 \s} \left( g^{mk} \wh E_k(\s) - \frac{\g^m}{\a} \ps_o(\s)\right) \ ,\\
\GGa {\qs_o}   {\qs_o} {\ps_o} & =\frac{1}{2 \a} \qs_o(\b)+   \frac{1}{2\a^2}\g^m \g^k  g_{mk} \qs_o(\a)-  \frac{1}{4\a^2}\g^m \g^k  g_{mk}\ps_o(\b)   -  \frac{1}{\a^2}\b\qs_o(\a)+ \nonumber\\
& \hskip 1 cm  +  \frac{\b}{2\a^2}\ps_o(\b) 
   -\frac{\g^m}{2\a} \qs_o(\g^i)g_{im} + \frac{\g^m}{4\a}\wh E_m(\b)-\nonumber \\
& \hskip 1 cm-  \frac{\b }{ 2\s}\left(\frac{1}{\a} \qs_o(\s) +  \frac{1}{2\a^2} \left(  \g^m \g^k  g_{mk} - 2\b \right) \ps_o(\s) - \frac{ \g^m}{2 \a} \wh E_m(\s)\right)   \ ,
\end{align}
\begin{align}
\label{772ter} \GGa {\qs_o}   {\qs_o} {\qs_o} & =  \frac{1}{\a}\qs_o(\a)-  \frac{1}{2\a}\ps_o(\b)+\frac{1}{\s} \qs_o(\s)- \frac{\b}{ 2\a \s}\ps_o(\s)\ .
  \end{align}
\end{theo}
\par 
\bigskip
The proof will be carried out in three steps, which we provide in the next subsections. In the first step we compute all covariant derivatives $\n_{X_A} X_B$   determined by  two vector fields of the adapted frame field $(X_A) = (\wh E_i, \ps_o, \qs_o)$ {\it under the assumption  $\s \equiv 1$}. In the second step, the  determined covariant derivatives are used   to compute  the  Christoffel symbols $\GGa A B C$, {\it still under the condition $\s \equiv 1$}. In the  concluding third step, the Christoffel symbols $\GGa A B C$ are determined with no restriction on $\s$ by  using   classical transformation formulas for the Levi-Civita covariant derivatives under conformal changes of the metric. 
\par
\medskip
\subsection{The first step}\label{casesigma1} 
By  Koszul's formula, for any triple of vector fields $X_1, X_2, X_3 $,  
\begin{multline} \label{Koszul}  g(\n_{X_1} X_2, X_3) = \frac{1}{2} \bigg( X_1(g(X_2, X_3)) + X_2(g( X_1, X_3)) - X_3(g(X_1, X_2)) - \\
  - g([X_1, X_3], X_2) - g(  [X_2, X_3], X_1) + g([X_1, X_2],X_3)\bigg) \ .\end{multline}
Using this formula, we may determine the functions $g(\n_{X_1} X_2, X_3)$,    for a compatible metric $g$ {\it with $\s \equiv 1$}, for  with any choice of  $X_1$, $X_2$, $X_3$  in   a set  of vector fields  of the form
$$\big\{\ \wh X, \ps_o, \qs_o, \ \text{where}\ \wh X\  \text{is the lift of a vector field}\ X\ \text{on}\ N\ \big\}\ .$$
We get the following expressions: 
 \begin{align} \label{primissima}
& \n_{\wh X} \wh Y:&  &  g(\n_{\wh X}  \wh Y, \wh Z) =  g_o(\n^o_{X}   Y, Z)  +  g(S_{XY}, Z)\ , \\
& &  & g(\n_{\wh X}  \wh Y, \ps_o) 
 = - \frac{\a}{4}  g_o(J X, Y)\ , \\
& &  & g(\n_{\wh X}  \wh Y, \qs_o) =   \frac{1}{4}  \wh X(\g^k g_o(Y, E_k)) + \frac{1}{4}  \wh Y(\g^k g_o( X,E_k) ) 
 -\frac{1}{4}  \b g_o(J X, Y)  \ ,
 \end{align}
 where $S$ is the tensor field of type $(0, 3)$  on $N$, defined by
 $$ g(S_{XY}, Z)  \= 
 \frac{\g^j}{4}g_o(J X, Z) g_o(E_j, Y) + \frac{\g^j}{4} g_o(J Y,  Z)  g_o(E_j, X) -  \frac{\g^j}{4} g_o(J X, Y)  g_o(E_j,Z)\ ;$$

 \begin{align}
&  \n_{\wh X} \ps_o:  &  & g(\n_{\wh X}  \ps_o, \wh Z) = 
   \frac{\a}{4}g_o(J X,  Z) \ ,\\
 && & g(\n_{\wh X}  \ps_o, \ps_o) =  0\ , \\
  && &  g(\n_{\wh X}  \ps_o, \qs_o) =   \frac{1}{4} \wh X(\a) +   \frac{1}{4}  \ps_o( \g^i) g_o(X, E_i)\ ; \hskip 4 cm 
 \end{align}
     \begin{align}
\nonumber  & \n_{\wh X} \qs_o :  & &g(\n_{\wh X}  \qs_o, \wh Z) = \frac{1}{4}  \wh X(\g^t g_o( E_t,  Z))  -\frac{1}{4}  \wh Z(\g^i g_o( X,  E_i)) -\\
 & & & \hskip 5cm  -\frac{1}{4}  \g^t g_o([X, Z], E_t)  +\frac{1}{4}  \b g_o(J X, Z)  \ ,  \hskip 0.5 cm 
  \\ 
 && &g(\n_{\wh X}  \qs_o, \ps_o) =  \frac{1}{4}  \wh X(\a)   -    \frac{1}{4} \ps_o(\g^i)  g_o(X, E_i) \ , \\
 && & g(\n_{\wh X}  \qs_o, \qs_o) =  \frac{1}{4} \wh X(\b)\ ;  \hskip 4 cm 
     \end{align}
 \begin{align}
 & \n_{\ps_o} \wh Y: && g(\n_{ {\ps_o}}  \wh Y, \wh Z) =    \frac{\a}{4}  g_o(J Y, Z)\ ,\\
 & &&  g(\n_{ {\ps_o}}  \wh Y, \ps_o) =  0\ ,\\
&& &  g(\n_{ {\ps_o}}  \wh Y, \qs_o) = \frac{1}{4}  \ps_o(\g^i) g_o(Y, E_i) +  \frac{1}{4}  \wh Y(\a) \ ;\hskip 4 cm 
 \end{align}
   \begin{align}
 & \n_{\ps_o} \ps_o: & &  g(\n_{ {\ps_o}}  \ps_o, \wh Z) =0 \ ,\\
 & & & g(\n_{ {\ps_o}}  \ps_o, \ps_o)  = 0\ ,\\
& &  &  g(\n_{ {\ps_o}}  \ps_o, \qs_o) =  \frac{1}{2} \ps_o(\a)\ ; \hskip 7 cm 
   \end{align}
   \begin{align}
 & \n_{\ps_o} \qs_o: & &  g(\n_{ {\ps_o}}  \qs_o, \wh Z) = \frac{1}{4}  \ps_o(\g^i)g_o(E_i, Z)  -  \frac{1}{4}  \wh Z(\a)\ ,\\
 &&  &  g(\n_{ {\ps_o}}  \qs_o, \ps_o) =  0\ ,\\
 &&  &  g(\n_{ {\ps_o}}  \qs_o, \qs_o) =  \frac{\ps_o(\b)}{4} \ ; \hskip 7 cm
   \end{align}
  \begin{align}
\nonumber  &\n_{\qs_o} \wh Y: &  &g(\n_{ {\qs_o}}  \wh Y, \wh Z) =  \frac{1}{4} \wh Y(\g^i g_o( E_i,  Z))   -  \frac{1}{4} \wh Z(\g^t g_o( Y,  E_t)) -\\
 & & & \hskip 5cm  - \frac{1}{4}  \g^t g_o([Y, Z], E_t) +  \frac{1}{4}  \b g_o(J Y, Z)  \ , \hskip 0.5 cm 
  \\ 
&&  & g(\n_{ {\qs_o}}  \wh Y, \ps_o) = \frac{1}{4} \big( \wh Y(\a) - \ps_o(\g^i) g_o(E_i, Y) \big)\ ,\\
&&  & g(\n_{ {\qs_o}}  \wh Y, \qs_o) = \frac{\wh Y(\b)}{4}   \ ;
  \end{align}
   \begin{align}
 & \n_{\qs_o} \ps_o: & &  g(\n_{ {\qs_o}}  \ps_o, \wh Z) = \frac{1}{4}  \ps_o(\g^i) g_o(E_i, \wh Z) -  \frac{1}{4} \wh Z(\a) \ ,\\
  && &  g(\n_{ {\qs_o}}  \ps_o, \ps_o) =  0\ ,\\
 &&  &   g(\n_{ {\qs_o}}  \ps_o, \qs_o) =   \frac{\ps_o (\b)}{4}\ ;\hskip 7 cm 
   \end{align}
    \begin{align}
 & \n_{\qs_o} \qs_o: &   & g(\n_{ {\qs_o}}  \qs_o, \wh Z) =  \frac{1}{2} \qs_o(\g^i) g_o(E_i, Z)  -  \frac{1}{4} \wh Z(\b)\ ,\\
 &&   & g(\n_{ {\qs_o}}  \qs_o, \ps_o) =  \frac{1}{2} \qs_o(\a) -  \frac{\ps_o(\b)}{4} \ ,\\
\label{ultimissima} &&  &  g(\n_{ {\qs_o}}  \qs_o, \qs_o) =   \frac{ \qs_o (\b)}{4} \ .\hskip 7 cm 
\end{align}
From this list, we may   recover the explicit expressions of the covariant derivatives of vector fields of the adapted frame field $(\wh E_i, \ps_o, \qs_o)$ as follows.  We  claim that
 the  dual coframe field $(\wh E^i, \ps^*_o, \qs^*_o)$ is given  by  the following $1$-forms (here,   $(g^{\ell m}) \= (g_{ij})^{-1} = \big( g_o(E_i, E_j)\big)^{-1}$)
 \begin{multline}\label{dualcoframe}  \wh E^i = g\bigg(g^{ik} \wh E_k - \frac{\g^i}{\a} \ps_o, \cdot\bigg) \ ,\quad
 \ps_o^* = g\bigg(\frac{2}{\a} \qs_o +  \frac{1}{\a^2} \left(  \g^m \g^k  g_{mk} - 2\b \right) \ps_o - \frac{ \g^m}{\a} \wh E_m, \cdot\bigg)\ ,\\
  \qs_o^* = g\bigg(\frac{2}{\a} \ps_o, \cdot\bigg)\ .
\end{multline}
 This claim  can be  checked  using \eqref{condiriso} and   observing that  the right hand sides in  the above equalities are $1$-forms that satisfy the equalities 
 \begin{align*}
 & \wh E^i(\wh E_j) = g^{ik} g_{kj}  = \d^i_j\ ,\quad  \wh E^i(\ps_o) = 0\ ,\quad   \wh E^i(\qs_o)  = g^{ik} \frac{\g^m}{2} g_{mk} -  \frac{\g^i}{\a} \frac{\a}{2} = 0\ ,\\
& \ps_o^*(\wh E_j) = \frac{2}{\a} \frac{\g^m}{2} g_{jm} - \frac{ \g^m}{\a} g_{mj} = 0\ ,\quad  \ps_o^*(\ps_o) = \frac{2}{\a} \frac{\a}{2} = 1\ ,\\   
&\hskip 4 cm \ps_o^*(\qs_o)  = \frac{2}{\a}\frac{\b}{2}+\frac{1}{\a^2} \left(  \g^m \g^k  g_{mk} - 2\b \right)  \frac{\a}{2} -  \frac{ \g^m}{\a} \frac{\g^k}{2} g_{mk} = 0\ ,\\
&  \qs_o^*(\wh E_j) =  0\ ,\quad  \qs_o^*(\ps_o) = 0\ ,\quad \qs_o^*(\qs_o)  = \frac{2}{\a} \frac{\a}{2} = 1\ .
\end{align*}
Since any local vector field $Z$ on $M$ can be written in terms of the frame field    $(\wh E_i, \ps_o, \qs_o)$ as 
$$Z = \wh E^i(Z) \wh E_i + \ps_o^*(Z) \ps_o + \qs_o^*(Z) \qs_o\ ,$$
from the above  expressions for the $1$-forms $\wh E^i, \ps^*_o$, and $\qs^*_o$,  
we get that for any pair of vector fields $X, Y$ on $M$, the Levi-Civita covariant derivative $\n_X Y$
is equal to  
\begin{multline}\label{3.48}  \n_X Y =   \left( g^{mk} g( \n_X Y ,  \wh E_k) -  \frac{\g^m}{\a} g(\n_X Y,  \ps_o) \right) \wh E_m +\\
+  \left(  \frac{2}{\a} g( \n_X Y, \qs_o) + \frac{1}{\a^2} \left(  \g^m \g^k  g_{mk} - 2\b \right)  g( \n_X Y, \ps_o) - \frac{ \g^m}{\a} g( \n_X Y, \wh E_m) \right)\ps_o + \\
+  \left(\frac{2}{\a}  g(\n_X Y, \ps_o)\right)\qs_o\ .
\end{multline}
Combining  \eqref{primissima} -- \eqref{ultimissima} with  \eqref{3.48}, we  get   the  covariant derivatives we are looking for.  We list   them in  \eqref{731} -- \eqref{739}   (here, we denote by $S_{ij|m}$ the components of the tensor field $S$ in terms of the frame field $(E_i)$ on $N$): 
 \begin{multline} \label{731}
\n_{\wh E_i} \wh E_j =  \left( 
  g^{mk} g_o(\n^o_{E_i}   E_j, E_k)  +  g^{mk}S_{ij|k}
  + \frac{\g^m \o_{ij} }{4}  \right) \wh E_m+\\
+  \left(    \frac{1}{2 \a} \wh E_i(\g^k g_{jk} )+   \frac{1}{2 \a}\wh E_j(\g^k g_{ik} ) - \xcancel{\frac{\b}{2\a } \o_{ij}} -  \right.\\
\left.  -  \frac{1}{4 \a} \g^m \g^k  g_{mk} \o_{ij} + \xcancel{\frac{2\b}{ \a^2} \frac{\a}{4} \o_{ij}}  - \frac{ \g^m}{\a} g_o(\n^o_{E_i}   E_j, E_m)  - \frac{ \g^m}{\a}  S_{ij|m} \right)\ps_o  -   \frac{  \o_{ij}}{2} \qs_o\ ,
 \end{multline}
\beq 
 \n_{\wh E_i} \ps_o =  \frac{\a g^{mk} \o_{ik} }{4}  \wh E_m
+  \left(    \frac{1}{2\a}  \wh E_i(\a) +  \frac{1}{2\a}  \ps_o( \g^k) g_{ik}  - \frac{ \g^m \o_{im}}{4}
 \right)\ps_o \ ,\hskip 2.3 cm 
\eeq
 \begin{multline}
 \n_{\wh E_i} \qs_o  = \Bigg(\frac{g^{mk}}{4} \wh E_i (\g^t g_{tk})- \frac{g^{mk}}{4}\wh E_k (\g^t g_{ti})  - \frac{\g^\ell}{4} c_{i r}^t g_{t \ell} g^{m r} +  \frac{g^{mk}}{4}\b \o_{ik}  -\\
 \hskip 7 cm - \frac{\g^m}{4\a}\wh E_i(\a)  +  \frac{\g^m}{4\a}\ps_o(\g^t)g_{it}\Bigg)\wh E_m+\\
 +\bigg( \frac{1}{2\a}\wh E_i(\b)+ \frac{1}{4\a^2} \g^m \g^k  g_{mk}  \wh E_i(\a)  -  \frac{1}{4\a^2} \g^m \g^k  g_{mk} \ps_o(\g^t)g_{it}  - \frac{\b}{2\a^2}\wh E_i(\a)  + \frac{\b}{2\a^2}\ps_o(\g^t)g_{it} -\\
  -\frac{\g^m}{4\a}\wh E_i(\g^tg_{tm}) + \frac{\g^m}{4\a}\wh E_m(\g^tg_{it})  +   \frac{\g^m\g^t}{4 \a} g_{t \ell} c^\ell_{i m} - \frac{\g^m}{4\a}\b\o_{im}\bigg)\ps_o +\\
   + \Bigg( \frac{1}{2\a}\wh E_i(\a) - \frac{1}{2\a}\ps_o(\g^t)g_{it}\Bigg)\qs_o\ ,
\end{multline}

\beq
  \n_{\ps_o} \wh E_i =  \frac{\a g^{mk}}{4}\o_{ik}\wh E_m+\Bigg( \frac{1}{2\a}\wh E_i(\a)+ \frac{1}{2\a}\ps_o(\g^t)g_{it}  -  \frac{\g^m}{4}\o_{im} \Bigg)\ps_o\ ,\hskip 2 cm 
 \eeq
\beq 
  \n_{\ps_o} \ps_o = 
 \ps_o(\log \a) \ps_o \ ,\hskip 9 cm 
\eeq
\begin{multline} 
  \n_{\ps_o} \qs_o =\bigg(  \frac{1}{4}\ps_o(\g^m) -  \frac{g^{mk}}{4}\wh E_k(\a)\Bigg)  \wh E_m+\Bigg( \frac{1}{2\a}\ps_o(\b) -\frac{\g^m}{4\a}\ps_o(\g^i)g_{im}+ \frac{\g^m}{4\a}\wh E_m(\a)   \bigg)\ps_o\ ,
\end{multline}
\begin{multline}
 \n_{\qs_o} \wh E_i= \bigg( \frac{g^{mk} }{4}\wh E_i(\g^tg_{tk})- \frac{g^{mk} }{4}\wh E_k(\g^tg_{ti}) - \frac{\g^\ell}{4} c_{i r}^t g_{t \ell} g^{m r}  +\frac{g^{mk} }{4}\b\o_{ik}  -\\
 \hskip 7 cm-\frac{\g^m}{4\a}\wh E_i(\a) + \frac{\g^m}{4\a}\ps_o(\g^t )g_{ti}\bigg)\wh E_m+\\
+ \bigg(\frac{1}{2\a}\wh E_i(\b)+  \frac{1}{4 \a^2} \g^m \g^k g_{mk} \wh E_i(\a)- \frac{1}{4 \a^2} \g^m \g^k g_{mk}\ps_o(\g^t )g_{it}  - \frac{1}{  2 \a^2} \b\wh E_i(\a) + \\
+ \frac{1}{   2\a^2} \b \ps_o(\g^t )g_{it} - \frac{\g^m}{4\a} \wh E_i(\g^tg_{tm}) +  \frac{\g^m}{4\a}\wh E_m(\g^tg_{it}) + \frac{\g^m \g^t}{4 \a} g_{t\ell} c^\ell_{i m}  - \frac{\g^m}{4\a}\b\o_{im}\bigg)\ps_o+\\
+\left(  \frac{1}{2 \a}  \wh E_i(\a)-  \frac{1}{2 \a} \ps_o(\g^t )  g_{it}\right)\qs_o\ ,
\end{multline}
\begin{multline}
  \n_{\qs_o} \ps_o =\Bigg(  \frac{g^{mk}}{4} \ps_o(\g^i)g_{ik}-  \frac{g^{mk}}{4}\wh E_k(\a)\Bigg)\wh E_m+\Bigg( \frac{1}{2\a}\ps_o(\b) -\frac{\g^m}{4\a}\ps_o(\g^t)g_{tm}+\frac{\g^m}{4\a}\wh E_m(\a) \Bigg)\ps_o\ ,
\end{multline}
\ \\[-20 pt]
\begin{multline} \label{739}
  \n_{\qs_o} \qs_o = \Bigg(\frac{g^{mk}}{2}\qs_o(\g^i)g_{ik} - \frac{g^{mk}}{4}\wh E_k(\b)  -\frac{\g^m}{2\a} \qs_o(\a) +\frac{\g^m}{4\a}\ps_o(\b)\Bigg)\wh E_m+\\
  +\Bigg(\frac{1}{ 2\a} \qs_o(\b)+   \frac{1}{2\a^2}\g^m \g^k  g_{mk} \qs_o(\a)-  \frac{1}{4\a^2}\g^m \g^k  g_{mk}\ps_o(\b)   -  \frac{\b}{\a^2}\qs_o(\a) +   \frac{\b}{2\a^2}\ps_o(\b) -\\
  -\frac{\g^m}{2\a} \qs_o(\g^i)g_{im} + \frac{\g^m}{4\a}\wh E_m(\b) \Bigg)\ps_o+ \left(  \frac{1}{\a}\qs_o(\a)-  \frac{\ps_o(\b)}{2\a}\right)\qs_o\ .
\end{multline}
\hfill
\par
\medskip
\subsection{The second step}
Let us now denote by $\Ga A B C$ the Christoffel symbols of the Levi-Civita connection of a compatible metric $g$ as in \eqref{buona}  under the assumption that the function $\s$ is identically equal to $1$.  Since the  $\Ga A B C$ are  the functions  that appear in  the expansions 
$ \n_{X_A} X_B = \Ga A B C X_C$ of  the covariant derivatives     \eqref{731} -- \eqref{739},  all such Christoffel symbols can be  determined by just looking at  those formulas.   For  convenience of the reader,   we  provide the complete  list in the next lines 
 \par
\begin{align}
\Ga i j m &=  
  g^{mk} g_o(\n^o_{E_i}   E_j, E_k)  +  g^{mk}S_{ij|k}
 + \frac{\g^m \o_{ij} }{4}  \ ,\\
\nonumber  \Ga i j {\ps_o} &=    \frac{1}{2 \a} \wh E_i(\g^k g_{jk} )+   \frac{1}{2 \a}\wh E_j(\g^k g_{ik} )  + \\
&\hskip 2 cm   -  \frac{1}{4 \a} \g^m \g^k  g_{mk} \o_{ij}   - \frac{ \g^m}{\a} g_o(\n^o_{E_i}   E_j, E_m)  - \frac{ \g^m}{\a}  S_{ij|m} \ ,\\
 \Ga i j {\qs_o}  &= - \frac{  \o_{ij}}{2} \ ,\\
  \Ga i {\ps_o} m & = \Ga {\ps_o} i m = \frac{\a g^{mk} \o_{ik} }{4} \ ,\\
 \Ga i {\ps_o} {\ps_o} &=   \Ga {\ps_o} i {\ps_o} =   \frac{1}{2\a}  \wh E_i(\a) +  \frac{1}{2\a}  \ps_o( \g^k) g_{ik}  - \frac{ \g^m \o_{im}}{4}\ ,
\\
 \Ga i {\ps_o} {\qs_o}& = \Ga {\ps_o} i {\qs_o} = 0\ , \\
 \nonumber \Ga i {\qs_o} m & =  \Ga  {\qs_o} i m  = \frac{g^{mk}}{4} \wh E_i (\g^t g_{tk})- \frac{g^{mk}}{4}\wh E_k (\g^t g_{ti})  - \frac{\g^\ell}{4} c_{i r}^t g_{t \ell} g^{m r}  +\frac{g^{mk}}{4}\b \o_{ik} -\\
& \hskip 7 cm - \frac{\g^m}{4\a}\wh E_i(\a)  + \frac{\g^m}{4\a}\ps_o(\g^t)g_{it} \ ,\end{align}
 \begin{align}
\nonumber  \Ga i {\qs_o} {\ps_o} & =  \Ga  {\qs_o} i  {\ps_o}  =  \frac{1}{2\a}\wh E_i(\b)+ \frac{1}{4\a^2} \g^m \g^k  g_{mk}  \wh E_i(\a) -  \frac{1}{4\a^2} \g^m \g^k  g_{mk} \ps_o(\g^t)g_{it}  -  \frac{\b}{2\a^2}\wh E_i(\a) +\\
 & \hskip 1 cm   +   \frac{1}{2\a^2}\b\ps_o(\g^t)g_{it}  -\frac{\g^m}{4\a}\wh E_i(\g^tg_{tm}) + \frac{\g^m}{4\a}\wh E_m(\g^tg_{it}) + \frac{\g^m \g^t}{4 \a} g_{t \ell} c^\ell_{i m} - \frac{\g^m}{4\a}\b\o_{im}\ ,
 \\
 \Ga i {\qs_o} {\qs_o} & =  \Ga  {\qs_o} i  {\qs_o}  = \frac{1}{2\a}\wh E_i(\a)  -  \frac{1}{2\a}\ps_o(\g^t)g_{it}\ ,\\
 \Ga {\ps_o} {\ps_o} m & = 0\ ,\\
  \Ga {\ps_o} {\ps_o} {\ps_o} & = \ps_o(\log \a) \ ,\\
  \Ga {\ps_o} {\ps_o} {\qs_o} &  = 0\ ,\\
   \Ga {\ps_o} {\qs_o} m & =   \Ga {\qs_o} {\ps_o} m =  \frac{1}{4}\ps_o(\g^m) -  \frac{g^{mk}}{4}\wh E_k(\a)\ ,\\
  \Ga {\ps_o} {\qs_o} {\ps_o} & =    \Ga {\qs_o} {\ps_o}  {\ps_o}  = \frac{1}{2\a}\ps_o(\b) -\frac{\g^m}{4\a}\ps_o(\g^i)g_{im}+ \frac{\g^m}{4\a}\wh E_m(\a)\ ,  \\
  \Ga {\ps_o} {\qs_o} {\qs_o} & =   \Ga {\qs_o}  {\ps_o} {\qs_o}  =  0\ ,\\[10pt]
  %
\Ga {\qs_o} {\qs_o} m & =  \frac{g^{mk}}{2}\qs_o(\g^i)g_{ik} - \frac{g^{mk}}{4}\wh E_k(\b)  -\frac{\g^m}{2\a} \qs_o(\a) +\frac{\g^m}{4\a}\ps_o(\b)\ ,
\\
\nonumber \Ga {\qs_o} {\qs_o} {\ps_o} & = \frac{1}{ 2\a} \qs_o(\b)+   \frac{1}{2\a^2}\g^m \g^k  g_{mk} \qs_o(\a)-  \frac{1}{4\a^2}\g^m \g^k  g_{mk}\ps_o(\b)   -  \frac{\b}{\a^2}\qs_o(\a) + 
\end{align}
\begin{align}
& \hskip 1 cm +  \frac{\b}{2\a^2} \ps_o(\b)
   -\frac{\g^m}{2\a} \qs_o(\g^i)g_{im} + \frac{\g^m}{4\a}\wh E_m(\b) \ m \ ,\\
 \Ga {\qs_o} {\qs_o} {\qs_o}   & =   \frac{1}{\a}\qs_o(\a)-  \frac{1}{2\a}\ps_o(\b)\ .
\end{align}
Note that  the equalities  $\Ga i {\ps_o}  A = \Ga {\ps_o} i A$, $\Ga i {\qs_o}  A = \Ga {\qs_o} i A$, etc. are also  consequences of the fact that the torsion of the Levi-Civita connection is $0$ and that the pairs of vector fields $\{\wh E_i$, $\ps_o\}$, $\{\wh  E_i$ $\qs_o\}$, etc.,   commute. \par
\medskip
\subsection{The third step}
Assume that   $g$ is one of the  metrics   considered in the previous two subsections (i.e., compatible with $\s \equiv 1$)  and denote by   $D$   the Levi-Civita connection of a conformally scaled   metric $g^\f= e^{2\f} g$ for some smooth  $\f$. It is  well known that,  for any pair of vector fields $X, Y$ of $M$ (see e.g. \cite{Be}*{Th. 1.159}),
\beq  \label{confD} D_X Y = \nabla_X Y + X(\f) Y + Y(\f) X -  g(X, Y)  \grad(\f) \ .\eeq 
If we expand  $\grad \f$  in terms of the frame field $(\wh E_i, \ps_o, \qs_o)$ as 
\beq \grad \f = (\grad \f)^{\wh E_i} \wh E_i +  (\grad \f)^{\ps_o} \ps_o + (\grad \f)^{\qs_o} \qs_o\ ,\eeq
we see that  the   Christoffel symbols $\Ga A B C$  for a compatible  metric $g$ with $\s \equiv 1$, as considered in the previous subsections,  and the Christoffel symbols $\GGa A B C$  for  the conformally scaled metric $g^\f$  are related to each other by 
\begin{align}
\label{746}
\GGa i j m &= \Ga i j m + \wh E_i(\f) \d_j^m + \wh E_j(\f) \d_i^m - g_{ij}  (\grad \f)^{\wh E_m}\ ,\\
\GGa i j {\ps_o} & =  \Ga i j {\ps_o}  - g_{ij}  (\grad \f)^{\ps_o}\ ,\\
\GGa i j {\qs_o} & =  \Ga i j {\qs_o} - g_{ij}  (\grad \f)^{\qs_o}\ ,
\\[7pt] 
\GGa i {\ps_o} m &= \GGa  {\ps_o} i m  = \Ga i {\ps_o} m + \ps_o(\f) \d_i^m \ ,
\\
\GGa i {\ps_o} {\ps_o} &=  \GGa  {\ps_o} i {\ps_o} = \Ga i {\ps_o} {\ps_o}  + \wh E_i(\f) \ ,\\
\GGa i {\ps_o} {\qs_o} &= \GGa  {\ps_o} i {\qs_o} = \Ga i {\ps_o} {\qs_o} \ , \\[7pt] 
\GGa i {\qs_o} m &= \GGa  {\qs_o} i m  =  \Ga i {\qs_o} m +\qs_o(\f)\d_i^m   -   \frac{\g^t}{2} g_{ti} (\grad \f)^{\wh E_m}\ ,
\\
\GGa i {\qs_o} {\ps_o} &=  \GGa  {\qs_o} i {\ps_o}  = \Ga i {\qs_o} {\ps_o}     -    \frac{\g^t}{2} g_{ti} (\grad \f)^{\ps_o} \ ,\\
\GGa i {\qs_o} {\qs_o} &= \GGa  {\qs_o} i {\qs_o}  =  \Ga i {\qs_o} {\qs_o}  + \wh E_i(\f)   -    \frac{\g^t}{2} g_{ti} (\grad \f)^{\qs_o}\ , \\[7pt] 
\GGa {\ps_o}  {\ps_o}  m & =  \Ga {\ps_o}  {\ps_o}  m  \ ,\\
\GGa {\ps_o}   {\ps_o} {\ps_o} & =  \Ga {\ps_o}  {\ps_o}  {\ps_o}+2\ps_o(\f)\ , \\
\GGa {\ps_o}   {\ps_o} {\qs_o} & =  \Ga {\ps_o}   {\ps_o} {\qs_o}\ , \end{align}
\begin{align}
\GGa {\ps_o}  {\qs_o}  m & = \GGa {\qs_o}   {\ps_o}   m  =  \Ga {\ps_o}  {\qs_o}  m-  \frac{\a }{2}(\grad \f)^{\wh E_m} \ , \\
\GGa {\ps_o}   {\qs_o} {\ps_o} & =  \GGa  {\qs_o}   {\ps_o} {\ps_o} = \Ga {\ps_o}   {\qs_o} {\ps_o}+\qs_o(\f)- \frac{\a}{2}(\grad \f)^{ {\ps_o}} \ , \\
\GGa {\ps_o}   {\qs_o} {\qs_o} & =  \GGa {\qs_o}  {\ps_o}   {\qs_o} = \Ga {\ps_o}   {\qs_o} {\qs_o} +\ps_o(\f)-  \frac{\a }{2}(\grad \f)^{ {\qs_o}} \ , \\[7pt] 
\GGa {\qs_o}  {\qs_o}  m & = \Ga {\qs_o}  {\qs_o}  m-   \frac{\b}{2}(\grad\f)^{\wh E_m} \ ,\\
\GGa {\qs_o}   {\qs_o} {\ps_o} & = \Ga {\qs_o}  {\qs_o}  {\ps_o}-    \frac{\b}{2}(\grad\f)^{\ps_o} \ ,  \\
\label{772} \GGa {\qs_o}   {\qs_o} {\qs_o} & =  \Ga {\qs_o}  {\qs_o}  {\qs_o}+2\qs_o(\f)-   \frac{\b}{2}(\grad\f)^{\qs_o} \ .
\end{align}
  We now recall that  any  vector field $X$ on $M$  decomposes into  the sum
\begin{equation}
\begin{split}  X & =   \wh E^i(X) \wh E_i + \ps_o^*(X) \ps_o + \qs_o^*(X) \qs_o  
 = g\left(X,  g^{ik} \wh E_k - \frac{\g^i}{\a} \ps_o\right) \wh E_i + \\
& \hskip 0.5 cm + g\left(X,\frac{2}{\a} \qs_o +  \frac{1}{\a^2} \left(  \g^m \g^k  g_{mk} - 2\b \right) \ps_o - \frac{ \g^m}{\a} \wh E_m\right) \ps_o 
+  g\left(X, \frac{2}{\a} \ps_o\right) \qs_o \ .\end{split}
\eeq
From this, we get that the components $(\grad \f)^A$ of the gradient of $\f$ are equal to 
\beq  \label{742*}
\begin{split}
& (\grad \f)^{\wh E_i} \= g^{ik} \wh E_k(\f) - \frac{\g^i}{\a} \ps_o(\f)\ ,\\
& (\grad \f)^{\ps_o} \=  \frac{2}{\a} \qs_o(\f) +  \frac{1}{\a^2} \left(  \g^m \g^k  g_{mk} - 2\b \right) \ps_o(f) - \frac{ \g^m}{\a} \wh E_m(f)\ ,\\
&  (\grad \f)^{\qs_o} \=  \frac{2}{\a} \ps_o(\f) .
\end{split}
\eeq
Inserting these expressions and \eqref{746} -- \eqref{772}  into \eqref{746} -- \eqref{772},  we get the explicit formulas for  the Christoffel symbols $\GGa A B C$ of the scaled metric $g^\f = e^{2 \f} g$. They are:
\begin{align}
\label{746bis}
\GGa i j m &=
  g^{mk} g_o(\n^o_{E_i}   E_j, E_k)  +  g^{mk}S_{ij|k}
  +\frac{\g^m \o_{ij} }{4}   + \wh E_i(\f) \d_j^m + \wh E_j(\f) \d_i^m\nonumber \\
 &\hskip 1 cm- g_{ij}  \left(  g^{mk} \wh E_k(\f) - \frac{\g^m}{\a} \ps_o(\f)\right)\ ,
\\
 \GGa i j {\ps_o}  &=  \frac{1}{2 \a} \wh E_i(\g^k g_{jk} )+  \frac{1}{2 \a}\wh E_j(\g^k g_{ik} )   
 -  \frac{1}{4 \a} \g^m \g^k  g_{mk} \o_{ij}   - \frac{ \g^m}{\a} g_o(\n^o_{E_i}   E_j, E_m)  - \frac{ \g^m}{\a}  S_{ij|m} \nonumber\\
 &\hskip 1 cm- g_{ij}  \left(\frac{2}{\a} \qs_o(\f) +  \frac{1}{\a^2} \left(  \g^m \g^k  g_{mk} - 2\b \right) \ps_o(\f) - \frac{ \g^m}{\a} \wh E_m(\f)\right)\ ,\\
 \GGa i j {\qs_o}  &=  -  \frac{  \o_{ij}}{2}  - \frac{2 g_{ij} }{\a} \ps_o(\f)\ ,\end{align}
\begin{align}
 \GGa i {\ps_o} m & =  \GGa  {\ps_o} i m =   \frac{\a g^{mk} \o_{ik} }{4} + \ps_o(\f) \d_i^m \ ,\\
 \GGa i {\ps_o} {\ps_o}&=  \GGa  {\ps_o} i {\ps_o} = \frac{1}{2\a}  \wh E_i(\a) +  \frac{1}{2\a}  \ps_o( \g^k) g_{ik}- \frac{ \g^m \o_{im}}{4} + \wh E_i(\f) \ ,\\
 \GGa i {\ps_o} {\qs_o}& =  \GGa i {\ps_o} {\qs_o} =0\ ,\\[7 pt]
 \GGa i {\qs_o} m &=  \GGa  {\qs_o} i m  = \frac{g^{mk} }{4}\wh E_i(\g^tg_{tk})- \frac{g^{mk} }{4}\wh E_k(\g^tg_{ti})    - \frac{\g^\ell}{4} c_{i r}^t g_{t \ell} g^{m r} +  \frac{g^{mk} }{4}\b\o_{ik} - \nonumber \\
&-\frac{\g^m}{4\a}\wh E_i(\a) + \frac{\g^m}{4\a}\ps_o(\g^t )g_{ti} +\qs_o(\f)\d_i^m-  \frac{\g^t }{2} g_{ti} \left(g^{mk} \wh E_k(\f) - \frac{\g^m}{\a} \ps_o(\f)\right)\ , \end{align}
 \begin{align}
 %
 \GGa i {\qs_o} {\ps_o} &=   \GGa {\qs_o} i {\ps_o} = \frac{1}{2\a}\wh E_i(\b)+  \frac{1}{4 \a^2} \g^m \g^k g_{mk} \wh E_i(\a) - \frac{1}{4 \a^2} \g^m \g^k g_{mk}\ps_o(\g^t )g_{it}  - \frac{1}{ 2\a^2} \b\wh E_i(\a) + \nonumber 
\\
& \hskip 1 cm +  \frac{1}{ 2\a^2} \b \ps_o(\g^t )g_{it} - \frac{\g^m}{4\a} \wh E_i(\g^tg_{tm}) +  \frac{\g^m}{4\a}\wh E_m(\g^tg_{it}) + \frac{\g^m \g^t}{4 \a} g_{t \ell} c^\ell_{i m} -  \frac{\g^m}{4\a}\b\o_{im} - \nonumber  \\
 & \hskip 1 cm-    \frac{\g^t}{2}  g_{ti} \left(\frac{2}{\a} \qs_o(\f) +  \frac{1}{\a^2} \left(  \g^m \g^k  g_{mk} - 2\b \right) \ps_o(\f) - \frac{ \g^m}{\a} \wh E_m(\f)\right)\ ,\\
 \GGa i {\qs_o} {\qs_o} &= \GGa  {\qs_o} i {\qs_o} =  \frac{1}{2\a}\wh E_i(\a)- \frac{1}{2\a}\ps_o(\g^t)g_{it}  + \wh E_i(\f)   -     \frac{\g^t}{2}  g_{ti} \left(\frac{2}{\a} \ps_o(\f)\right) \ , 
  \end{align}
 \begin{align} %
 %
 \GGa {\ps_o}  {\ps_o}  m & =  0\ , \\
 \GGa {\ps_o}   {\ps_o} {\ps_o} & =  \ps_o(\log \a)+2\ps_o(\f)\ ,\\
 \GGa {\ps_o}   {\ps_o} {\qs_o} & =  0 \ , \\ %
 \GGa {\ps_o}  {\qs_o}  m & =  \GGa  {\qs_o}  {\ps_o}   m  = \frac{1}{4}\ps_o(\g^m) -  \frac{g^{mk}}{4}\wh E_k(\a)-  \frac{\a }{2} \left (g^{mk} \wh E_k(\f) - \frac{\g^m}{\a} \ps_o(\f) \right)\ ,\\
 \GGa {\ps_o}   {\qs_o} {\ps_o} & =  \GGa  {\qs_o}  {\ps_o}  {\ps_o} = \frac{1}{2\a}\ps_o(\b) -\frac{\g^m}{4\a}\ps_o(\g^i)g_{im}+ \frac{\g^m}{4\a}\wh E_m(\a)+\qs_o(\f)-\nonumber\\
 & \hskip 1 cm-  \frac{\a}{2} \left(\frac{2}{\a} \qs_o(\f) +  \frac{1}{\a^2} \left(  \g^m \g^k  g_{mk} - 2\b \right) \ps_o(\f) - \frac{ \g^m}{\a} \wh E_m(\f) \right) \ ,  
 \\
  \GGa {\ps_o}   {\qs_o} {\qs_o} & = \GGa {\qs_o}   {\ps_o}   {\qs_o}  =\xcancel{\ps_o(\f)}- \xcancel{\frac{\a}{2}  \left ( \frac{2}{\a} \ps_o(\f)\right)}  = 0\ ,\\[7 pt]
\GGa {\qs_o}  {\qs_o}  m & = \frac{g^{mk}}{2}\qs_o(\g^i)g_{ik} - \frac{g^{mk}}{4}\wh E_k(\b)  -\frac{\g^m}{2\a} \qs_o(\a) +\frac{\g^m}{4\a}\ps_o(\b)-\nonumber\\
& \hskip 1 cm-  \frac{\b}{2}\left( g^{mk} \wh E_k(\f) - \frac{\g^m}{\a} \ps_o(\f)\right) \ , \\
\GGa {\qs_o}   {\qs_o} {\ps_o} & = \frac{1}{ 2\a} \qs_o(\b)+   \frac{1}{2\a^2}\g^m \g^k  g_{mk} \qs_o(\a)-  \frac{1}{4\a^2}\g^m \g^k  g_{mk}\ps_o(\b)   -  \frac{1}{\a^2}\b\qs_o(\a) + \nonumber\end{align}
\begin{align}
& \hskip 1 cm +  \frac{1}{2\a^2}\b \ps_o(\b) 
   -\frac{\g^m}{2\a} \qs_o(\g^i)g_{im} + \frac{\g^m}{4\a}\wh E_m(\b)-\nonumber \\
& \hskip 1 cm-  \frac{\b}{2}\left(\frac{2}{\a} \qs_o(\f) +  \frac{1}{\a^2} \left(  \g^m \g^k  g_{mk} - 2\b \right) \ps_o(\f) - \frac{ \g^m}{\a} \wh E_m(\f)\right)  \ , \\
\label{772bis} \GGa {\qs_o}   {\qs_o} {\qs_o} & =  \frac{1}{\a}\qs_o(\a)-  \frac{\ps_o(\b)}{2\a}+2\qs_o(\f)- \frac{\b}{2}  \left(\frac{2}{\a} \ps_o(\f)\right)\ .
  \end{align}

\par
In order to conclude, it is now sufficient to  observe that the metric \eqref{buona} with an arbitrary  $\sigma> 0$  can be obtained from the metric considered in  \S \ref{casesigma1} (i.e., with $\s \equiv 1$)  by  applying the scaling factor $e^{2\f}$ with  $\f \= \frac{1}{2} \log \s$.  Hence, the desired expressions for the Christoffel symbols are  given by  \eqref{746bis} -- \eqref{772bis} with $\f$ replaced by $ \frac{1}{2} \log \s$ at all places.   These substitutions yield  \eqref{746ter} -- \eqref{772ter}.
\par

\par
\newpage

\vskip 1.5truecm
\hbox{\parindent=0pt\parskip=0pt
\vbox{\baselineskip 9.5 pt \hsize=3.5truein
\obeylines
{\smallsmc
Dmitri V.  Alekseevsky
Institute for Information Transmission Problems
B. Karetny per. 19
127051 Moscow
Russia
\&
University of Hradec  Kr\'alov\'e,
Faculty of Science, 
Rokitansk\'eho 62, 
500~03 Hradec Kr\'alov\'e,
Czech Republic
%
%
%
\
}\medskip
{\smallit E-mail}\/: {\smalltt dalekseevsky@iitp.ru}
}
\vbox{\baselineskip 9.5 pt \hsize=3.1truein
\obeylines
{\smallsmc
Masoud Ganji
School of Science and Technology
University of New England,
Armidale NSW 2351
Australia
\phantom{\&}
\phantom{University of Hradec   Kr{\'a}lov\'e,}
\phantom{Faculty of Science, Rokitansk\'eho 62, 500~03 Hradec Kr\'alov\'e,}
\phantom{Czech Republic}
\phantom{Czech Republic}
\phantom{Czech Republic}
\phantom{Czech Republic}}\medskip
{\smallit E-mail}\/: {\smalltt mganjia2@une.edu.au}
}
}
\vskip 1truecm
\hbox{\parindent=0pt\parskip=0pt
\vbox{\baselineskip 9.5 pt \hsize=3.5truein
\obeylines
{\smallsmc
Gerd Schmalz
School of Science and Technology
University of New England,
Armidale NSW 2351
Australia
\
}\medskip
{\smallit E-mail}\/: {\smalltt schmalz@une.edu.au}
}
\vbox{\baselineskip 9.5 pt \hsize=3.1truein
\obeylines
{\smallsmc
Andrea Spiro
Scuola di Scienze e Tecnologie
Universit\`a di Camerino
Via Madonna delle Carceri
I-62032 Camerino (Macerata)
Italy
}\medskip
{\smallit E-mail}\/: {\smalltt andrea.spiro@unicam.it
}
}
}

\end{document}